\documentclass[12pt]{iopart}          
\begin{document}
\title
[A gradient system on the quantum information space
realizing the Karmarkar flow]
{A gradient system on the quantum information space that realizes
the Karmarkar flow for linear programming }
\author{
Yoshio Uwano
and
Hiromi Yuya
}
\address{
Department of Complex Systems, Future University-Hakodate,
\\
116-2, Kameda Nakano-cho, Hakodate, 041-8655, Japan
}
\ead{uwano@fun.ac.jp}
\begin{abstract}
In the paper of Uwano
[Czech. J. of Phys., vol.56, pp.1311-1316 (2006)],
a gradient system is found on the space of density matrices endowed
with the quantum SLD Fisher metric (to be referred to as the quantum
information space) that realizes a generalization of a gradient system
on the space of multinomial distributions studied by Nakamura
[Japan J. Indust. Appl. Math., vol.10, pp.179-189 (1993)].
On motived by those papers, the present paper aims to construct
a gradient system on the quantum information space that realizes the
Karmarkar flow, the continuous limit of the Karmarkar projective
scaling algorithm for linear programming.
\end{abstract}
%
%
%
\maketitle
\newcommand{\dbc}[2]{\langle {#1} , {#2} \rangle}
\newcommand{\Mnm}{\mbox{M}(2^n ,m)}
\newcommand{\Mnmone}{\mbox{M}_{1}(2^n ,m)}
\newcommand{\trace}{\mbox{{\rm tr}} \;}
\newcommand{\Herpos}{\dot{P}_m}
\newcommand{\Herposemi}{P_m}
\newcommand{\Rmet}[2]{(\!( {#1} , {#2} )\!)^{R}}
\newcommand{\QFmet}[2]{(\!( {#1} , {#2} )\!)^{QF}}
\newcommand{\CFmet}[2]{(\!( {#1} , {#2} )\!)^{CF}}
\newcommand{\Smet}[2]{(\!( {#1} , {#2} )\!)^{Smp}}
\newcommand{\Dmet}[2]{(\!( {#1} , {#2} )\!)^{D}}
\newcommand{\inv}{\pi_m^{-1}(\Herpos)}
\newcommand{\regS}{\dot{{\cal S}}}
\newcommand{\pim}{\pi_{m}}
\newcommand{\rank}{\mbox{{\rm rank}} \,}
\newcommand{\bsl}{\backslash}
\newtheorem{definition}{Definition}[section]
\newtheorem{theorem}[definition]{Theorem}
\newtheorem{lemma}[definition]{Lemma}
\newtheorem{conjecture}[definition]{Conjecture}
\newtheorem{remark}{Remark}
\newtheorem{proposition}[definition]{Proposition}
\newtheorem{corollary}[definition]{Corollary}
\section{Introduction}
There exist various excellent algorithms developed in engineering
and systems science; as one of those, the Karmarkar projective scaling
algorithm for linear programming is very famous, which solves linear
programming problems in polynomial time (Karmarkar 1984).
\par
In the middle of 1990's, Nakamura revealed an integrability behind
several algorithms and governing equations arising in engineering
and systems science; the Karmarkar flow emerging as a continuous
version of the Karmarkar projective scaling algorithm (Karmarkar 1990,
Nakamura 1994a), an algorithm solving the eigenvalue problem
of anti-Hermitean matrices (Nakamura 1992), an averaged learning
equation of the Hebb-type (Nakamura 1994b) and a gradient system
relevant to the information geometry of multinomial distributions
(Nakamura 1993).
\par
More than ten years after Nakamura's work (Nakamura 1993) on the
gradient system on the space of multinomial distributions (GS-MD),
one of the authors, Y.~Uwano, encountered a very natural generalization
of the GS-MD through a study on geometry and dynamics of a search
algorithm for an ordered-tuple of multi-qubit states (Uwano 2006,
Uwano \etal 2007):
The space of density matrices endowed with the quantum SLD Fisher
metric which is referred to as the quantum information space (QIS)
was constructed through a geometric reduction of the space of
ordered-tuples of multi-qubit states. On the QIS, the gradient
system associated with the negative von Neumann entropy (GS-NVNE)
was studied: The matrix-form solution of the GS-NVNE was obtained
explicitly, whose diagonal part was shown to describe the solution
of the GS-MD. The GS-NVNE on the QIS is therefore understood as the
natural generalization of the GS-MD.
\par
This encounter encourages the authors to seek other gradient systems
on the QIS that realize certain dynamical systems or algorithms in
engineering and systems science.
The aim of the present paper is to construct a gradient system on
the QIS that realizes the Karmarkar flow for linear programming.
In what follows, the organization of this paper is outlined.
\par
In section 2, the general framework of gradient systems on the QIS
is derived together with a review of the QIS: The construction of the
QIS is accomplished by applying a geometric reduction method to the
space of ordered tuples of multi-qubit states for a quantum search
(Uwano 2006, Uwano \etal 2007). After the review, the equation of motion
is derived for arbitrary gradient systems on the QIS.
Section~3 is the core part of the present paper, where the
gradient system on the QIS realizing the Karmarkar flow (GS-QIS-K)
is given explicitly. The Karmarkar flow mentioned here is the
family of trajectories arising from the Karmarkar projective scaling
algorithm for the {\it unconstrained case} (Karmarkar 1990, Nakamura
1993). A key to find the GS-QIS-K is to observe that
the Riemannian structure of the canonical simplex, the underlying
manifold, for the Karmarkar flow is in isometry to the Riemannian
submanifold of the QIS consisting of diagonal matrices.
Section~4 is devoted to the concluding remarks.
Throughout the present paper, differential geometric calculus works
very effectively, for whose detail Appendices A and B are prepared.    
\section{Gradient systems on the QIS}
\subsection{Geometric setting-up for the QIS: Review}
Following Uwano (2006) and Uwano \etal (2007), we review the geometric
reduction method to construct the quantum information space (QIS)
from {\it the space of normalized orederd-tuples of multi-qubit states}
(STMQ). The natural Riemannian metric of the STMQ is shown to be reduced
to the quantum SLD-Fisher metric of QIS.
Those who are not so familiar to differential geometry may skip
subsection~2.1, which is however deserves for tracing basic
ideas for a series of works on gradient systems on the QIS by the
author(s) (see also the closing remark in subsection~2.1).
\subsubsection{Reduction of the STQM to the QIS}
\quad
By $\Mnm$, we denote the Hilbert space of $2^n\times m$ complex
matrices endowed with a natural Hermitean inner product
\begin{eqnarray}
\label{def-inner}
\langle \Phi, \Phi^{\prime} \rangle 
=
\frac{1}{m} \tr(\Phi^{\dagger}\Phi^\prime)
\quad (\Phi, \Phi^{\prime} \in \Mnm),
\end{eqnarray}
where the superscript ${}^\dagger$ will indicate the Hermitean conjugate
operation from now on. The $\Mnm$ is thought to describe the Hilbert
space of ordered tuples of multi-qubit states, if every column,
$\phi^{(j)} \in {\bf C}^{2^n}$ ($j=1,\,2,\,\ldots,\,m$), of a
matrix 
\begin{equation}
\label{Phi}
\Phi=\bigl(\phi^{(1)},\,\ldots,\,\phi^{(m)}\bigr)\in\Mnm
\end{equation}
is understood to express a vector in the standard complex Hilbert space,
${\bf C}^{2^n}$, of $n$-qubit states (see Nielsen and Chuang 2000
for basic setting-up of quantum computation, e.g., multi-qubit states
\ldots {\it etc}).
\begin{remark}
The set of $q\times r$ complex matrices is denoted by $M(q,r)$
through this paper.
\end{remark}
We denote by $\Mnmone$ the subset of $\Mnm$ consisting of $\Phi$'s in
$\Mnm$ with the unit norm. Namely,
\begin{eqnarray}
\label{def-Mnmone}
\Mnmone = \{
\Phi \in \Mnm \, \vert \, \langle \Phi , \Phi \rangle =1 \} .
\end{eqnarray}
The space $\Mnmone$ is understood as that of normalized ordered-tuples
of multi-qubit states, which will be abbreviated to the STMQ.
\par
Let us consider the space of $m \times m$ density matrices
(see Nielsen and Chuang 2000),
\begin{eqnarray}
\label{P_m}
\Herposemi = 
\{ \rho \in M(m,m) \, \vert \, 
\rho^{\dag} = \rho , \, \tr \rho =1 , \,
\rho : \mbox{positive semidefinite} \} ,
\end{eqnarray}
as the quotient space, say $\Mnmone / \mathrm{U}(2^n)$, of $\Mnmone$
with respect to a natural left $\mathrm{U}(2^n)$ action
\begin{eqnarray}
\label{U(2^n)}
\Phi \in \Mnmone \longmapsto g \Phi \in \Mnmone
\quad
(\Phi \in \Mnmone, \, g \in \mathrm{U}(2^n)),
\end{eqnarray}
where $\mathrm{U}(2^n)$ denotes the group of unitary matrices of
degree $2^n$.
Indeed, if we define the map of $\Mnmone$ to $\Herposemi$ to be
\begin{eqnarray}
\label{pi_m}
\pi_m : \Phi \in \Mnmone \longmapsto
\frac{1}{m} \Phi^{\dag}\Phi \in \Herposemi ,
\end{eqnarray}
we see that $\pi_m$ is surjective
and that $\pi_m (\Phi)=\pi_m (\Phi^{\prime})$ holds true
for $\Phi, \Phi^{\prime} \in \Mnmone$ if and only if there
exists $g \in \mathrm{U}(2^n)$ subject to $\Phi=g\Phi^{\prime}$.
This shows $\pi_m (\Mnmone)= \Herposemi \cong \Mnmone /\mathrm{U}(2^n)$.
It follows from (\ref{Phi}) and (\ref{U(2^n)}) that the
$\mathrm{U}(2^n)$ action leaves {\it relative-configuration}
among $\phi^{(j)}$s (multi-qubit states) in an ordered tuple
$\Phi$ invariant, so that $\Herposemi$ realizing
$\Mnmone/\mathrm{U}(2^n)$ can be referred to as {\it the space of
relative-configurations of multi-qubit states in ordered tuples}
(Uwano 2006 and Uwano \etal 2007).
\subsubsection{The quantum SLD Fisher metric on the QIS}
\quad
To proceed differential calculus including differential
equations, we have to consider a regular part of $\Herposemi$,
which is realized as the $m \times m$ regular density matrices
denoted by $\Herpos$ (Uwano 2006 and Uwano \etal 2007), which can be
referred to as
{\it the space of regular relative-configurations of multi-qubit
states in ordered tuples}.
\par
As the space of $m \times m$ regular density matrices, 
$\Herpos$ admits the quantum information space structure endowed
with the quantum SLD Fisher metric denoted by $\QFmet{\cdot}{\cdot}$.
The quantum SLD Fisher metric is introduced as follows.
\par
Let us consider the space of $m \times m$ traceless Hermitean matrices,
\begin{equation}
\label{tan-P_m} T_{\rho}\Herpos=\bigl\{\Xi\in M(m,m)\,\vert\,
\Xi^{\dagger}=\Xi,\,\tr\Xi=0\bigr\}\,,
\end{equation}
as the tangent space of $\Herpos$ at $\rho$. The symmetric logarithmic
derivative (SLD) at $\rho\in\Herpos$ for $\Xi\in T_{\rho}\Herpos$ is
defined to provide the matrix $\mathcal{L}_{\rho}(\Xi)\in
M(m,m)$ subject to
\begin{equation}
\label{SLD} 
\frac{1}{2}\,\bigl\{\rho\mathcal{L}_{\rho}(\Xi)+
\mathcal{L}_{\rho}(\Xi)\,\rho\bigr\}=\Xi\qquad (\Xi\in
T_{\rho}\Herpos)\, .
\end{equation}
The quantum SLD Fisher metric, denoted by $\QFmet{\cdot}{\cdot}$,
is then defined to be
\begin{equation}
\label{defeq-Fisher}
\fl
\quad
\QFmet{\Xi}{\Xi^{\prime}}_{\rho}=\frac{1}{2}
\tr \left[\rho\bigl(L_{\rho}(\Xi)L_{\rho}(\Xi^{\prime})+
L_{\rho}(\Xi^{\prime})L_{\rho}(\Xi)\bigr)\right]
\qquad
(\Xi,\Xi^{\prime}\in T_{\rho}\Herpos) \, ,
\end{equation}
(see Amari and Nagaoka 2000, Uwano \etal 2007).
\par
We wish to express $\QFmet{\cdot}{\cdot}$ explicitly.
Let $\rho \in \Herpos$ be expressed as
\begin{eqnarray}
\begin{array}{l}
\rho=h\Theta h^{\dagger}, \quad h \in\mathrm{U}(m)
\\ \noalign{\medskip}
\label{rho}
\Theta=\mathrm{diag}(\theta_1,\ldots,\theta_m)
\quad
\mbox{with}
\quad
\tr \Theta =1 , \quad \theta_{k}>0 \,\, (k=1,2,\cdots m) ,
\end{array}
\end{eqnarray}
where $\mathrm{U}(m)$ denotes the group of $m \times m$ unitary
matrices. Expressing $\Xi \in T_{\rho} \Herpos$ as
\begin{eqnarray}
\label{Xi-chi}
\Xi = h \chi h^{\dagger}
\end{eqnarray}
with $h \in \mathrm{U}(m)$ in (\ref{rho}), we obtain an
explicit expression,
\begin{eqnarray}
\label{SLD-exp}
(h^{\dagger} {\cal L}_{\rho}(\Xi) h)_{jk}
= 
\frac{2}{\theta_j + \theta_k} \chi_{jk} \quad (j,k=1,2,\cdots m),
\end{eqnarray}
of the SLD to $\Xi \in T_{\rho}\Herpos$.
Putting (\ref{rho})-(\ref{SLD-exp}) into (\ref{defeq-Fisher}), we have
\begin{equation}
\label{exp-QF}
\QFmet{\Xi}{\Xi^{\prime}}_{\rho}
=
2\sum_{j,k=1}^{m}
\frac{\overline{\chi}_{jk}\chi^{\prime}_{jk}}{\theta_j+\theta_k}
\end{equation}
where $\Xi^{\prime} \in T_{\rho}\Herpos$ is expressed as
\begin{eqnarray}
\label{Xi'-chi'}
\Xi^{\prime} = h \chi^{\prime} h^{\dagger}.
\end{eqnarray}
\par
Surprisingly, the quantum SLD Fisher metric $\QFmet{\cdot}{\cdot}$ thus
defined turns out to be identical with the Riemannian metric, denoted by
$\Rmet{\cdot}{\cdot}$, that is determined along with the dimensional
reduction of $\inv \, (\subset \Mnmone)$ to $\Herpos$ through $\pi_m$
given by (\ref{pi_m}).
Indeed, under (\ref{rho}), (\ref{Xi-chi}) and (\ref{Xi'-chi'}),
we have
\begin{eqnarray}
\label{exp-Rmet}
\Rmet{\Xi}{\Xi^{\prime}}_{\rho}
=
\frac{1}{2} \sum_{j,k=1}^{m}
\frac{\overline{\chi}_{jk}\chi^{\prime}_{jk}}{\theta_j+\theta_k}
=
\frac{1}{4} \QFmet{\Xi}{\Xi^{\prime}}_{\rho}
\end{eqnarray}
(see Appendix~A for the construction of $\Rmet{\cdot}{\cdot}$
and Uwano \etal 2007). 
\begin{theorem}
[Uwano 2006, Uwano \etal 2007]
The quantum SLD Fisher metric $\QFmet{\cdot}{\cdot}$ coincides
with $\Rmet{\cdot}{\cdot}$ up to the constant multiple;
\begin{eqnarray}
\label{QF-R}
\QFmet{\Xi}{\Xi^{\prime}}_{\rho}
=
4 \, \Rmet{\Xi}{\Xi^{\prime}}_{\rho}
\quad
(\Xi , \Xi^{\prime} \in T_{\rho}\Herpos ).
\end{eqnarray}
\end{theorem}
\smallskip
Throughout this paper, we will refer to the space of $m \times m$
regular density matrices, $\Herpos$, endowed with the quantum SLD Fisher
metric $\QFmet{\cdot}{\cdot}$ as the quantum information space (QIS)
for brevity, that will be indicated also as the pair
$(\Herpos, \QFmet{\cdot}{\cdot})$. 
\subsection{General framework of gradient systems on the QIS}
Now that the QIS is constructed as the Riemannian manifold with
the metric $\QFmet{\cdot}{\cdot}$, we move on to study gradient
systems on the QIS. The motive already mentioned in Sec.~1 for
studying gradient systems on the QIS comes from the
big similarity discovered by Uwano (2006) between the gradient system
on the QIS with the potential equal to the negative von Neumann entropy
and that on the information space of multinomial distributions studied by
Nakamura (1993).
\par
Let $F$ be a smooth real-valued function on $\Herpos$.
With the quantum SLD Fisher metric $\QFmet{\cdot}{\cdot}$,
the gradient vector, denoted by $(\mathrm{grad}\, F)(\rho)$,
of $F$ at $\rho$ is defined to be the tangent vector of $\Herpos$
at $\rho$ subject to
\begin{eqnarray}
\label{grad}
\QFmet{(\mathrm{grad} \, F)(\rho)}{\Xi^{\prime}}_{\rho}
=
(\mathrm{d}F)_{\rho}(\Xi^{\prime})
=
\left. \frac{d}{dt} \right\vert_{t=0} F(\gamma(t))
\quad ({}^{\forall}\Xi^{\prime} \in T_{\rho}\Herpos ),
\end{eqnarray}
(see Kobayashi and Nomizu 1969). The $\gamma (t)$ with
$a \leq t \leq b$ ($a<0<b$) in (\ref{grad}) is a curve in
$\Herpos$ subject to
\begin{eqnarray}
\label{curve-gamma}
t \in [a, \, b] \mapsto \gamma (t) \in \Herpos ,
\quad
\gamma (0)=\rho ,
\quad
\left. \frac{d\gamma}{dt} \right\vert_{t=0}= \Xi^{\prime} .
\end{eqnarray}
Through a straight but quite long calculation, an explicit expression
of the gradient vector $(\mathrm{grad} \, F)(\rho)$ is obtained in
what follows.
\par
Let us introduce the real-valued variables,
$x_{jk}, \,\, y_{jk}$ ($1 \leq j<k \leq m$) and
$z_{\ell}$ ($\ell =1,2,\cdots,m$),
to express the entries of $\rho \in \Herpos$ as
\begin{eqnarray}
\label{rho-xyz}
\fl
\rho_{jk}=\overline{\rho}_{kj}= x_{jk}+iy_{jk}
\quad (1 \leq j<k \leq m), \quad 
\rho_{\ell \ell}=z_{\ell} \quad (\ell =1,2,\cdots,m),
\end{eqnarray}
where $z_{\ell}$s are subject to
\begin{eqnarray}
\label{z-cond}
z_{\ell} >0 \quad (\ell =1,2,\cdots,m)
\quad \mbox{and} \quad
\sum_{\ell=1}^{m} z_{\ell}=1 .
\end{eqnarray}
For the complex variables $\rho_{ab}$ $(1 \leq a<b \leq m)$,
we introduce the partial differentiations,
\begin{eqnarray}
\label{diff-rho}
\fl
\frac{\partial}{\partial \rho_{ab}}
=
\frac{1}{2}\Big( 
\frac{\partial}{\partial x_{ab}} -i \frac{\partial}{\partial y_{ab}}
\Big) ,
\quad
\frac{\partial}{\partial {\overline \rho}_{ab}}
=
\frac{1}{2}\Big( 
\frac{\partial}{\partial x_{ab}} + i \frac{\partial}{\partial y_{ab}}
\Big)
 \quad (1 \leq a<b \leq m).
\end{eqnarray}
To calculate the rhs of (\ref{grad}), it is very convenient
to define the matrix-valued operator ${\cal M}$ to $F$ by
\begin{eqnarray}
\label{cal-M}
\big( {\cal M}(F) \big)_{jk}
=
\left\{ \begin{array}{ll}
\displaystyle{
\frac{\partial F}{\partial \overline{\rho}_{jk}}
=
\overline{ \frac{\partial F}{\partial \rho_{jk}}}
}
& \quad (1 \leq j < k \leq m)
\\
\noalign{\smallskip}
\displaystyle{
\frac{\partial F}{\partial \rho_{kj}}
}
& \quad (1 \leq k < j \leq m)
\\
\noalign{\smallskip}
\displaystyle{
\frac{\partial F}{\partial \rho_{jj}}
}
& \quad (j=k=1,2, \cdots ,m) .
\end{array}
\right.
\end{eqnarray}
By (\ref{cal-M}), the rhs of (\ref{grad}) is calculated to be
\begin{eqnarray}
\fl
\nonumber
\left. \frac{d}{dt} \right\vert_{t=0} F(\gamma(t))
=
\sum_{1 \leq a<b \leq m}
\left\{ \frac{\partial F}{\partial x_{ab}}( \rho )
        \Re (\Xi^{\prime}_{ab})
        +
        \frac{\partial F}{\partial y_{ab}}( \rho )
        \Im (\Xi^{\prime}_{ab})
\right\}
+
\sum_{r=1}^{m}
\frac{\partial F}{\partial z_{r}}( \rho ) \, \Xi^{\prime}_{rr}
\\
\nonumber
=
\sum_{1 \leq a<b \leq m}
\left\{ \frac{\partial F}{\partial x_{ab}}( \rho ) 
  \frac{1}{2} \Big( \Xi^{\prime}_{ab} + \overline{\Xi^{\prime}_{ab}} \Big)
+
\frac{\partial F}{\partial y_{ab}}( \rho ) 
  \frac{1}{2i} \Big( \Xi^{\prime}_{ab} - \overline{\Xi^{\prime}_{ab}} \Big)
\right\}
\\
\nonumber
\qquad
+
\sum_{r=1}^{m}
\frac{\partial F}{\partial z_{r}}(\rho) \, \Xi^{\prime}_{rr}
\\
\nonumber
=
\sum_{1 \leq a<b \leq m}
\left\{ \frac{\partial F}{\partial \rho_{ab}}(\rho) \Xi^{\prime}_{ab}
+
\frac{\partial F}{\partial \overline{\rho}_{ab}}(\rho) 
\overline{\Xi^{\prime}_{ab}} 
\right\}
+
\sum_{r=1}^{m}
\frac{\partial F}{\partial z_{r}}(\rho) \, \Xi^{\prime}_{rr}
\\
\nonumber
=
\sum_{1 \leq a<b \leq m}
\left\{ ({\cal M}(F))_{ba} \Xi^{\prime}_{ab}
+
({\cal M}(F))_{ab} \Xi^{\prime}_{ba}
\right\}
+
\sum_{r=1}^{m}
({\cal M}(F))_{rr} \Xi^{\prime}_{rr}
\\
\label{dF}
=
\tr \Big( {\cal M}(F) \Xi^{\prime} ) ,
\end{eqnarray}
where the symbols, $\Re$ and $\Im$ , stand for the real part and
the imaginary part, respectively.
\par
To calculate the lhs of (\ref{grad}),
the introduction of the Hermitean matrix,
\begin{eqnarray}
\label{cal-G}
{\cal G}=h^{\dagger} \, \big( (\mathrm{grad}F)(\rho) \big) \, h,
\end{eqnarray}
is of great use. Indeed, Eq.~(\ref{cal-G}) is put together with 
(\ref{rho}), (\ref{exp-QF}) and (\ref{Xi'-chi'}) to show
\begin{eqnarray}
\nonumber
\fl
\QFmet{(\mathrm{grad} \, F)(\rho)}{\Xi^{\prime}}_{\rho}
=
2 \sum_{j,k=1}^{m} 
\frac{\overline{{\cal G}}_{jk} \chi^{\prime}_{jk}}{\theta_j + \theta_k}
=
2 \sum_{j,k=1}^{m} 
\frac{{\cal G}_{kj}}{\theta_j + \theta_k}
\Big( \sum_{a,b=1}^{m} \overline{h}_{aj} \Xi^{\prime}_{ab} h_{bk} \Big)
\\
\label{lhs-grad}
=
2 \sum_{a,b=1}^{m} \Big( 
\sum_{j,k=1}^{m} h_{bk} {\tilde {\cal G}}_{kj} (h^{\dagger})_{ja} \Big)
\Xi^{\prime}_{ab}
=
2 \tr \Big( (h {\tilde {\cal G}} h^{\dagger}) \Xi^{\prime} \Big) ,
\end{eqnarray}
where ${\tilde {\cal G}}=({\tilde {\cal G}}_{jk})$ is the Hermitean matrix
defined to be
\begin{eqnarray}
\label{tilde-cal-G}
{\tilde {\cal G}}_{jk} = \frac{{\cal G}_{jk}}{\theta_j + \theta_k}
\quad (j,k=1,2,\cdots ,m).
\end{eqnarray}
Since the Hermitean form,
\begin{eqnarray}
\label{Herm-inner}
(\Xi, \Xi^{\prime}) \in T_{\rho} \Herpos \times T_{\rho} \Herpos
\mapsto
\tr (\Xi^{\dagger} \Xi^{\prime} ) \in {\bf C},
\end{eqnarray}
is well-known to be non-degenerate, we easily see, from
(\ref{dF}) and (\ref{lhs-grad}), that the equation
\begin{eqnarray}
\label{M-G}
{\cal M}(F)=2h {\tilde {\cal G}} h^{\dagger} + 2 cI 
\end{eqnarray}
has to hold true, where $c$ in is the constant emerging
from the trace-vanishing property of $\Xi^{\prime} \in T_{\rho}\Herpos$
(see (\ref{tan-P_m})). The value of $c$ will be determined below, soon.
From (\ref{cal-G}) and (\ref{M-G}), ${\cal G}$ turns out to take the form
\begin{eqnarray}
\nonumber
{\cal G}_{jk}
&=&
\frac{1}{2} (\theta_j + \theta_k)
\left\{ (h^{\dagger} {\cal M}(F) h )_{jk} -2c \delta_{jk} \right\}
\\
\label{G_jk}
&=& \frac{1}{2} 
\left( 
\Theta h^{\dagger} {\cal M}(F) h  
+
h^{\dagger} {\cal M}(F) h \Theta - 2c \Theta \right)_{jk}
\quad (j,k=1,2, \cdots),
\end{eqnarray}
where $\delta_{jk}$ denotes the Kronecker delta.
Equation (\ref{G_jk}) is combined with (\ref{cal-G}) and (\ref{rho})
to show
\begin{eqnarray}
\label{grad-F}
(\mathrm{grad} \, F)(\rho)
=
\frac{1}{2} \Big( 
\rho {\cal M}(F) + {\cal M}(F) \rho \Big) - c \rho.
\end{eqnarray}
We are now in a position to evaluate the constant $c$ by taking the
trace in both sides of (\ref{grad-F}). By a simple calculation,
the $c$ is determined to be
\begin{eqnarray}
\label{c}
c= \tr \Big( \rho {\cal M}(F) \Big),
\end{eqnarray}
that leads us to
\begin{eqnarray}
\label{grad-F-final}
(\mathrm{grad} \, F)(\rho)
=
\frac{1}{2} \Big( 
\rho {\cal M}(F) + {\cal M}(F) \rho \Big) 
- \Big( \tr \big( \rho {\cal M}(F) \big) \Big) \rho.
\end{eqnarray}
To summarize, we have the following.
\begin{lemma}
\label{lemma-grad}
If a gradient system on the quantum information space
$(\Herpos , \QFmet{\cdot}{\cdot})$ is associated with a potential
function $F$, it is governed by the equation of motion
\begin{eqnarray}
\label{grad-eq-F}
\frac{d \rho}{dt}
=
- \frac{1}{2} \Big( 
\rho {\cal M}(F) + {\cal M}(F) \rho \Big) 
+ \Big( \tr \big(\rho {\cal M}(F) \big) \Big) \rho,
\end{eqnarray}
where ${\cal M}(F)$ is the matrix defined by (\ref{cal-M}).
\end{lemma}
\begin{remark}
In the case of the GS-NVNE for example, Lemma~2.2 is not so effective
to derive its gradient vector since the negative von Neumann entropy
is hardly written in terms of the entries of $\rho$.
Indeed, in Uwano \etal 2007, the gradient vector was not calculated
directly on the QIS but was done through its {\it lifting} to STQM and
{\it projecting} to the QIS; For such calculation, the geometric
setting reviewed plays a central role. This could provide a good
account for our review in subsection 2.1.
\end{remark}
\section{The Karmarkar flow in the QIS}
This is the core section of this paper, where the Karmarkar
flow is realized on the QIS.
\subsection{The Karmarkar flow: Review}
Following Nakamura (1994), we review the Karmarkar flow for the
canonical linear programming problem
\begin{eqnarray}
\nonumber
\mathrm{minimize} & \quad & \sum_{j=1}^{m} c_j x_j
\\
\label{can-LP}
\mbox{subject to} & \quad & \sum_{j=1}^{m} x_j =1, \,
\quad  x_j \geq 0 \, \, \,  (j=1,2,\cdots ,m)
\end{eqnarray}
of {\it unconstrained case} (Karmarkar 1990), where $c_j$s are given
nonvanishing constants. We note here that no additional linear constraint
specifying a subspace of the $m-1$ dimensional canonical simplex,
\begin{eqnarray}
\label{simplex}
{\cal S}= \Big\{ x \in {\bf R}^m \, \Big\vert \,
\sum_{j=1}^{m}x_j =1 , \, x_j \geq 0 \, (j=1,2,\cdots,m)  \Big\} ,
\end{eqnarray}
is taken into account in the unconstrained case.
\par
A continuous limit of the Karmarkar projective scaling algorithm 
gives rise to the system of differential equations
\begin{eqnarray}
\label{DE-Karmarkar}
\frac{dx_j}{dt}=
-c_jx_j^2 + x_j \Big( \sum_{k=1}^{m}c_k x_k^2 \Big) \quad (j=1,2,\cdots,m).
\end{eqnarray}
The family of trajectories governed by (\ref{DE-Karmarkar}) is
what we refer to as the Karmarkar flow in the present paper.
The system of differential equations (\ref{DE-Karmarkar})
is brought into the following gradient-system form, according to
Nakamura (1993): To be free from singularities, our discussion will
be made on the regular part,
\begin{eqnarray}
\label{regS}
\regS = \Big\{ x \in {\bf R}^m \, \Big\vert \,
\sum_{j=1}^{m}x_j =1, \, x_j > 0 \, (j=1,2,\cdots,m) \Big\} ,
\end{eqnarray}
of the simplex ${\cal S}$ henceforth. With $\regS$, we endow the
Riemannian metric
\begin{eqnarray}
\label{Smet}
\Smet{u}{u^{\prime}}_{x} = \sum_{j=1}^{m} \frac{u_j u^{\prime}_j}{x_j}
\quad (u, u^{\prime} \in T_x \regS)
\end{eqnarray}
where the tangent space, denoted by $T_{x} \regS$, of $\regS$ at
$x \in \regS$ is given by
\begin{eqnarray}
\label{tan-regS}
T_{x}\regS
=
\Big\{ u \in {\bf R}^m \, \Big\vert \,
\sum_{j=1}^{m} u_j =0 \Big\}.
\end{eqnarray}
On $\regS$, the system of differential equations (\ref{DE-Karmarkar})
admits the gradient-system form if we take the function,
\begin{eqnarray}
\label{pot-Karmarkar}
\kappa (x) = \frac{1}{2} x^T C x \quad (x \in \regS )
\end{eqnarray}
as the potential, where $C$ is the diagonal matrix of the form
\begin{eqnarray}
\label{C}
C= \mathrm{diag} \, (c_1 , c_2 , \cdots , c_m) .
\end{eqnarray}
The gradient vector $(\mathrm{grad} \, \kappa)(x)$ at
$x \in \regS$ is defined to satisfy
\begin{eqnarray}
\label{def-grad-kappa}
\Smet{(\mathrm{grad} \, \kappa)(x)}{u^{\prime}}_{x}
=
\left. \frac{d}{dt} \right\vert_{t=0} \kappa (\sigma (t))
\quad ({}^{\forall} u^{\prime} \in T_{x} \regS ),
\end{eqnarray}
where $\sigma (t)$ with $a \leq t \leq b$ ($a<0<b$) is a curve in
$\regS$ subject to
\begin{eqnarray}
\label{curve-sigma}
t \in [a, \, b] \mapsto \sigma (t) \in \regS ,
\quad
\sigma (0)= x ,
\quad
\left. \frac{d\sigma}{dt} \right\vert_{t=0}= u^{\prime} 
\end{eqnarray}
(cf. (\ref{grad}) with (\ref{curve-gamma})).
By a straightforward calculation analogous to that for
$(\mathrm{grad}\, F)(\rho)$ in Sec.~2, we obtain
\begin{eqnarray}
\label{grad-kappa}
\big( (\mathrm{grad} \, \kappa)(x) \big)_j
=
c_jx_j^2 - x_j \Big( \sum_{k=1}^{m}c_k x_k^2 \Big) \quad (j=1,2,\cdots,m),
\end{eqnarray}
which yields (\ref{DE-Karmarkar}). See Appendix~B for a detail of the
calculation.
\subsection{The gradient system on the QIS realizing the Karmarkar flow}
The gradient system on the QIS that we are seeking is constructed
in what follows.
\subsubsection{The Riemannian structures of $\regS$ and the QIS}
Let us consider the submanifold,
\begin{eqnarray}
\label{cal-D}
\fl
{\cal D}
=
\Big\{ \rho \in \Herpos \,\, \Big\vert \,\,
\rho = \mathrm{diag} \, (\theta_{1}, \cdots , \theta_{m}) , \,
\sum_{k=1}^{m}\theta_k =1, \,\, \theta_{k} >0 \,\,\, (k=1,2,\cdots ,m)
\Big\} ,
\end{eqnarray}
of the QIS, which is easily seen to be diffeomorphic to $\regS$,
the regular part of the canonical simplex ${\cal S}$. Indeed, we
can find the smooth one-to-one and onto map,
\begin{eqnarray}
\label{mu}
\mu : x \in \regS
\mapsto \mathrm{diag} \, (x_{1}, \cdots , x_{m}) \in {\cal D}
\subset \Herpos .
\end{eqnarray}
Restricting the quantum SLD-Fisher metric $\QFmet{\cdot}{\cdot}$
of the QIS to the submanifold ${\cal D}$, we can make ${\cal D}$
the Riemannian submanifold, whose metric will be denoted
by $\Dmet{\cdot}{\cdot}$ henceforth. Namely, on expressing
the tangent space of ${\cal D}$ at $\Theta$ (cf. (\ref{rho}))
as the subspace,
\begin{eqnarray}
\label{tan-cal-D}
T_{\Theta}{\cal D}
=
\Big\{ Z \in M(m,m) \, \Big\vert \, 
Z =  \mathrm{diag} \, (\zeta_{1}, \cdots , \zeta_{m}),
\,
\sum_{j=1}^{m} \zeta_j = 0 \Big\},
\end{eqnarray}
of $T_{\Theta}\Herpos$, $\Dmet{\cdot}{\cdot}$ is defined to satisfy 
\begin{eqnarray}
\label{QF-D}
\Dmet{Z}{Z^{\prime}}_{\Theta}
=
\QFmet{Z}{Z^{\prime}}_{\Theta} \quad
(Z,Z^{\prime} \in T_{\Theta}{\cal D} \subset T_{\Theta}\Herpos).
\end{eqnarray}
We show the following.
\begin{lemma}
\label{regS-D}
The map $\mu$, defined by (\ref{mu}), of $\regS$ to ${\cal D}$
is isometric; the identity,
\begin{eqnarray}
\label{isometry}
\Dmet{\mu_{\ast ,x}(u)}{\mu_{\ast ,x}(u^{\prime})}_{\mu(x)}
=
\Smet{u}{u^{\prime}}_{x}
\quad (u, u^{\prime} \in T_{x}\regS ),
\end{eqnarray}
holds true, where $\mu_{\ast ,x}$ is the differential of the map
$\mu$ at $x$ defined by
\begin{eqnarray}
\label{mu_*x}
\mu_{\ast, x} (u^{\prime})
=
\left. \frac{d}{dt} \right\vert_{t=0} \mu (\sigma(t))
=
\mathrm{diag} \, (u^{\prime}_{1}, \cdots , u^{\prime}_{m})
\quad (u^{\prime} \in T_{x} \regS )
\end{eqnarray}
with (\ref{tan-regS}) and (\ref{curve-sigma}).
\end{lemma}
{\it Proof:}\quad
Equation (\ref{QF-D}) is put together with
Eqs.~(\ref{rho})-(\ref{Xi'-chi'}) and (\ref{mu_*x}) to yield
\begin{eqnarray}
\fl
\nonumber
\Dmet{\mu_{\ast ,x}(u)}{\mu_{\ast ,x}(u^{\prime})}_{\mu(x)}
=
\QFmet
{\mu_{\ast, x}(u)}{\mu_{\ast ,x}(u^{\prime})}_{\mu (x)}
\\ \noalign{\smallskip}
\label{calc-Dmet-Smet}
=
2\sum_{j,k=1}^{m}
\frac{\overline{(\mu_{\ast,x}(u))}_{jk}
      \, (\mu_{\ast,x}(u^{\prime}))_{jk}}{x_j+x_k}
=
\sum_{j=1}^{m}
\frac{u_{j} u^{\prime}_{j}}{x_j}
=
\Smet{u}{u^{\prime}}_{x}.
\end{eqnarray}
This completes the proof.
\subsubsection{Construction of the gradient system}
From Lemma~\ref{regS-D}, we learn the coincidence
of the Riemannian structures of the regular part, $\regS$, of the
canonical simplex for the Karmarkar flow, and of the submanifold,
${\cal D}$, of the QIS. Accordingly, in order to find a gradient
system realizing the Karmarkar flow on the QIS, we naturally come to
seek a function $K(\rho)$ on the QIS whose restriction to ${\cal D}$
coincides with the potential $\kappa (x)$ for the Karmarkar flow
through the map $\mu$:
We take $K(\rho)$ to be
\begin{eqnarray}
\label{K}
K(\rho)= \frac{1}{2} \tr \Big( C\rho^2 \Big),
\end{eqnarray}
where $C$ is the diagonal matrix given in (\ref{C}).
The entries of the matrix ${\cal M}(K)$ given by (\ref{cal-M})
with $K$ in place of $F$ are calculated to be
\begin{eqnarray}
\fl
\label{dK/drho_jk}
\left( {\cal M}(K) \right)_{jk}
=
\overline{ \frac{\partial K}{\partial \rho_{jk}} }
=
\overline{ \frac{1}{2} (c_j + c_k ) \rho_{kj} }
=
\frac{1}{2} (c_j + c_k ) \rho_{jk}
\quad (j,k=1, \cdots , m),
\end{eqnarray}
which is brought into the form
\begin{eqnarray}
\label{dK/drho}
\frac{\partial K}{\partial \rho}
=
\frac{1}{2} \big( C \rho + \rho C \big).
\end{eqnarray}
Equation (\ref{dK/drho}) is put together with (\ref{grad-F-final})
to show
\begin{eqnarray}
\label{grad-K}
(\mathrm{grad} \, K)(\rho)
=
\frac{1}{4}(\rho^2 C + 2 \rho C \rho + C \rho^2 )
- \Big( \tr \big( \rho C \rho) \Big) \rho ,
\end{eqnarray}
so that we have the following.
\begin{lemma}
\label{lemma-grad-K}
The gradient system on the QIS associated with the potential $K(\rho)$
is governed by the equation of motion,
\begin{eqnarray}
\label{grad-eq-K}
\frac{d \rho}{dt}
=
-\frac{1}{4}(\rho^2 C + 2 \rho C \rho + C \rho^2 )
+ \Big( \tr \big( \rho C \rho) \Big) \rho .
\end{eqnarray}
\end{lemma}
We are at the final stage to see how the gradient system on the QIS
with $K$ realize the Karmarkar flow. Recalling Eq.~(\ref{grad-K}),
we immediately obtain
\begin{eqnarray}
\label{grad-K-D}
(\mathrm{grad} \, K)(\Theta)
=
C \Theta^2 - \Big( \tr \big( C \Theta^2 ) \Big) \Theta
\in T_{\Theta} {\cal D} \subset T_{\Theta} \Herpos
\end{eqnarray}
(see (\ref{tan-cal-D}) for $T_{\Theta} {\cal D}$), which enables us
to restrict the equation of motion (\ref{grad-eq-K}) to the submanifold
${\cal D}$ of the QIS isometric to $\regS$. The restriction indeed gives
rise to the system of differential equations,
\begin{eqnarray}
\label{grad-eq-K-D}
\frac{d \theta_j}{dt}
=
-c_j \theta_j^2 + \theta_j \Big( \sum_{k=1}^{m} c_k \theta_k^2 \Big)
\quad
(j=1, 2, \cdots ,m),
\end{eqnarray}
on ${\cal D}$, which is evidently identical with the Karmarkar flow
(\ref{DE-Karmarkar}) with $\theta$ in place of $x$. In conclusion,
we have the following.
\begin{theorem}
\label{GS-K-Karmarkar}
The gradient system on the QIS associated with the potential
$K(\rho)$ (GS-QIS-K) realizes the Karmarkar flow on the submanifold
${\cal D}$ of the QIS.
\end{theorem}
\section{Concluding remarks}
We have successfully constructed the gradient system
(GS-QIS-K) on the QIS which realizes the Karmarkar flow on the
submanifold ${\cal D}$. A key to the success is the isometry of
the underlying Riemannian manifold $\regS$ for the Karmarkar flow
and the Riemannian submanifold ${\cal D}$ of the QIS, that is
presented in Lemma~\ref{regS-D}. 
\par
Through the study leading us to Lemma~\ref{regS-D} on the Riemannian
structures of $\regS$ and ${\cal D}$, a clear account for the encounter
with the GS-NVNE as a generalization of the GS-MD is obtained:
The Riemannian metric and the potential for the GS-MD are, respectively,
equal to those for the GS-NVNE restricted on ${\cal D}$ up to a common
multiplier.
\par
Integrability of the GS-QIS-K is an open question.
We wish to recall that, in the case of the GS-NVNE (Uwano 2006),
the invariance of the negative von Neumann entropy chosen as the potential
under the $\mathrm{U}(m)$ action,
$\rho \mapsto h \rho h^{\dagger}$ ($h \in \mathrm{U}(m)$), works
effectively to show the integraility in the sense that the GS-NVNE
allows a sufficient number of mutually independent integrals of motion.
The $\mathrm{U}(m)$ invariance of the potential $K (\rho)$
does not hold true, however, so that we have little expectation of the
integrability for the GS-QIS-K. If we find a $\mathrm{U}(m)$-invariant
potential whose restriction to ${\cal D}$ realize $\kappa (x)$, the
gradient system with that potential on the QIS would be integrable
and, further, would admit a double-Lax bracket representation like
in the case of the GS-NVNE.
\par
In the case that the matrix $C$ in the potential $K(\rho)$ is
taken to be $C=2I$, the GS-QIS-K has two other particular features:
One is that it turns out to realize not only the Karmarkar flow but
also the flow solving the eigenvalue problem of anti-Hermitean matrices
in view of Nakamura (1992, 1993). The other is that the potential $K(\rho)
=\tr (\rho^2)$ is called {\it the purity} whose logarithm provides the
minus of the quantum Renyi potential $ \log ( \tr \rho^q )  / (1-q)$
with $q=2$; the larger the purity of a quantum state is,
the larger its Hilbert-Schmidt distance from the maximally mixed state.
(Bengtsson and {\. Z}yczkowski 2006). Due to the second
feature, the potential $K(\rho)$ with $C=2I$ could be
interpreted to be an object in quantum physics. A paper dealing with
the case of $C=2I$ is in preparation.
\par
The averaged learning equation of Hebb-type dealt with in
Nakamura (1994) is a current target of the authors: Through another
geometric trick, we have recently succeeded to find its generalization
on the QIS, which will be reported in other paper (in preparation).
\par\noindent
\ack
\quad
The authors thank Dr.~Fumitaka Yura at Future University Hakodate for his
valuable comment on the physical meanings of the potential $K(\rho)$
with $C=2I$.
\appendix
\setcounter{section}{0}
\section{The Riemannian metric $\Rmet{\cdot}{\cdot}$}
\subsection{Geometry of $\inv$} 
We start with introducing the natural Riemannian metric of $\Mnmone$
($\supset \inv$).
On regarding $\Mnm \cong {\bf C}^{2^n}$ as the $2^{n+1}$-dimensional
Euclidean space, the real part of the Hermitean inner product of $\Mnm$
given by (\ref{def-inner}) can provide the Euclidean metric
\begin{eqnarray}
\nonumber
&&
(X, X^{\prime})^E_{\Phi} 
= 
\frac{1}{2m} \tr 
\big( X^{\dag}X^{\prime} + \overline{(X^{\dag} X^{\prime})} \big)
\\
\label{Euc-met}
&&
\qquad \qquad \qquad
(X, X^{\prime} \in \Mnm \cong T_{\Phi}\Mnm,
\, \Phi \in \Mnm) ,
\end{eqnarray}
where $T_{\Phi}{\bf C}^{2^n}$ denotes the tangent space of
$\Mnm$ at $\Phi \in \Mnm$. 
The tangent space, $T_{\Phi}\Mnmone$, of $\Mnmone$ at $\Phi \in \Mnmone$
is thereby defined to be
\begin{eqnarray}
\label{tan-Mnmone}
T_{\Phi}\Mnmone 
=
\{ X \in \Mnm \, \vert \, \tr (X^{\dag}\Phi + \Phi^{\dag}X) = 0 \},
\end{eqnarray}
on looking $\Mnmone$ upon as a submanifold of $\Mnm$.
Through the restriction of $T_{\Phi}\Mnm$ to $T_{\Phi}\Mnmone$
($\Phi \in \Mnmone$), the Euclidean metric $(\cdot , \cdot)^E$ of
$\Mnm$ is brought to the Riemannian metric of of $\Mnmone$,
\begin{eqnarray}
\label{R-met-Mnmone}
\fl
(X, X^{\prime})^R_{\Phi} 
= 
\frac{1}{2m} \tr
\big( X^{\dag}X^{\prime} + \overline{(X^{\dag} X^{\prime})} \big)
\qquad
(X, X^{\prime} \in T_{\Phi}\Mnmone, \, \Phi \in \Mnmone) .
\end{eqnarray}
On account that $T_{\Phi} \inv = T_{\Phi} \Mnmone$ for
$\Phi \in \inv$ ($\subset \Mnmone$), the $(\cdot , \cdot)^R$ becomes
the Rimennian metric of the inverse image $\inv$ of
$\Herpos$ by $\pi_m$ if restricted.
\subsection{The horizontal lift}
We introduce the horizontal lift of tangent vectors of
$T_{\rho}\Herpos$ to $T_{\Phi} \inv$ ($\pi_m (\Phi) = \rho \in \Herpos$)
as follows:
According to the $\mathrm{U}(2^n)$ action (\ref{U(2^n)}),
let us consider the orthogonal direct-sum decomposition
\begin{eqnarray}
\label{tan-decomp}
T_{\Phi}\Mnmone = \mathrm{Ver}(\Phi) \oplus_{\perp} \mathrm{Hor}(\Phi)
\quad (\Phi \in \inv)
\end{eqnarray}
for the metric $(\cdot , \cdot)^R_{\Phi}$ with
\begin{eqnarray}
\label{Ver}
\fl
\mathrm{Ver}(\Phi)
=
\{X \in T_{\Phi} \inv \, \vert \, X= \eta \Phi , \, \eta \in u(2^n) \}
\quad (\Phi \in \inv)
\end{eqnarray}
and
\begin{eqnarray}
\label{Hor}
\fl
\mathrm{Hor}(\Phi)
=
\{X \in T_{\Phi} \inv \, \vert \,
\Phi X^{\dag} - X \Phi^{\dag} = O_{2^n,2^n} \} 
\quad (\Phi \in \inv),
\end{eqnarray}
where $u(2^n)$ denotes the set of all the anti-Hermitean matrices
of degree-$2^n$ and $O_{2^n,2^n}$ does the null matrix of degree-$2^n$
(see also Prop.~4 in Uwano \etal 2007). 
The horizontal lift, denoted by $\ell_{\Phi}(\Xi)$, of
$\Xi \in T_{\rho}\Herpos$ is then defined to be
the unique tangent vector at $\Phi \in \pi_m^{-1}(\rho)$ subject to
\begin{eqnarray}
\label{lift-def}
{\pi_m}_{*, \Phi}(\ell_{\Phi}(\Xi))=\Xi 
\quad \mbox{and} \quad
\ell_{\Phi}(\Xi) \in \mathrm{Hor}(\Phi) .
\end{eqnarray}
The ${\pi_m}_{*, \Phi}$ is the differential of the map
$\pi_m$ at $\Phi \in \Herpos$, which is defined to be
\begin{eqnarray}
\label{pi_m*}
{\pi_m}_{*, \Phi} (X)
=
\left. \frac{d}{dt} \right\vert_{t=0} \pi_m (c(t))
\quad (X \in T_{\Phi} \Herpos),
\end{eqnarray}
where $c(t)$ with $a \leq t \leq b$ ($a<0<b$) is a curve in $\Mnmone$
satisfying
\begin{eqnarray}
\label{curve-c}
t \in [a, \, b] \mapsto c(t) \in \Mnmone,
\quad
c(0)=\Phi,
\quad
\left. \frac{dc}{dt} \right\vert_{t=0}=X .
\end{eqnarray}
Under the expressions (\ref{rho}) and (\ref{Xi-chi}) of
$\rho \in \Herpos$ and $\Xi \in T_{\rho}\Herpos$, the singular-value
decomposition (see Rao and Mitra 1971),
\begin{eqnarray}
\label{inv-rho}
\fl
\Phi 
=
g \, \left(
\begin{array}{c} \sqrt{m} \, \sqrt{\Theta} \\ O_{2^n-m,m} \end{array}
\right)\, h^{\dag}
\quad ({}^{\exists}g \in \mathrm{U}(2^n))
\quad \mbox{with} \quad
\sqrt{\Theta}=\mathrm{diag} (\sqrt{\theta_1}, \cdots , \sqrt{\theta_m}),
\end{eqnarray}
of $\Phi \in \pi_m^{-1}(\rho)$ works effectively to obtain the
horizontal lift.
Indeed, on putting (\ref{inv-rho}) with (\ref{rho}) and (\ref{Xi-chi})
together, the horizontal lift of $\Xi \in T_{\rho}\Herpos$ is given by
\begin{eqnarray}
\label{hor-lift-exp}
\ell_{\Phi}(\Xi)
=
\frac{\sqrt{m}}{2} g \, \left(
\begin{array}{c}
(\sqrt{\Theta})^{-1} ( \chi + \alpha_{\sqrt{\Theta}}( \chi))
\\
O_{2^n-m,m}
\end{array}
\right)
h^{\dag}
\end{eqnarray}
where $\alpha_{\sqrt{\Theta}}(\chi)$ stands for the $m \times m$
anti-Hermitean matrices uniquely determined by
\begin{eqnarray}
\label{alpha}
\Theta^{-1} \alpha_{\sqrt{\Theta}}(\chi)
+
\alpha_{\sqrt{\Theta}}(\chi) \Theta^{-1}
=
-\Theta^{-1}\chi + \chi \Theta^{-1}
\end{eqnarray}
(see Uwano \etal 2007). By a straightforward calculation, we have
\begin{eqnarray}
\label{exp-alpha}
\left( \alpha_{\sqrt{\Theta}}(\chi) \right)_{jk}
=
\left( \frac{\theta_{j}-\theta_{k}}{\theta_{j}+\theta_{k}} \right) \chi_{jk}
\quad (j,k = 1,2, \cdots , m).
\end{eqnarray}
\subsection{Defining $\Rmet{\cdot}{\cdot}$}
On using the Riemannian metric $(\cdot , \cdot)^R$ of
$\Mnmone$ and the horizontal lift $\ell_{\Phi}(\cdot)$,
the Riemannian metric $\Rmet{\cdot}{\cdot}$ of $\Herpos$
is defined to be
\begin{eqnarray}
\label{def-Rmet}
\Rmet{\Xi}{\Xi^{\prime}}_{\rho}
=
(\ell_{\Phi}(\Xi), \ell_{\Phi}(\Xi^{\prime}))_{\Phi}^R 
\quad (\rho \in \Herpos, \, \Xi , \Xi^{\prime} \in T_{\rho}\Herpos )
\end{eqnarray}
where $\Phi \in \pi_m^{-1}(\rho)$ can be chose arbitrarily.
Equations (\ref{R-met-Mnmone}), (\ref{hor-lift-exp}), (\ref{exp-alpha})
and (\ref{def-Rmet}) are put together to yield (\ref{exp-Rmet}).
\section{The gradient vectors}
We calculate the gradient vector $(\mathrm{grad} \, \kappa)(x)$
at $x \in \regS$ of $\kappa$ according to (\ref{def-grad-kappa}).
On taking (\ref{curve-sigma}) for the curve $\sigma (t)$ into account,
the rhs of (\ref{def-grad-kappa}) is calculated to be
\begin{eqnarray}
\label{rhs-def-grad-kappa}
\left. \frac{d}{dt} \right\vert_{t=0} \kappa (\sigma (t))
=
\sum_{j=1}^{m}
 \frac{\partial \kappa}{\partial x_j}(\sigma (0))
 \frac{d \sigma_j}{dt}(0)
=
\sum_{j=1}^{m}
 \frac{\partial \kappa}{\partial x_j}(x) u^{\prime}_j .
\end{eqnarray}
Further, Eq.~(\ref{Smet}) is combined with the lhs of
(\ref{def-grad-kappa}) to show
\begin{eqnarray}
\label{lhs-def-grad-kappa}
\Smet{(\mathrm{grad} \, \kappa)(x)}{u^{\prime}}_x
=
\sum_{j=1}^{m}
 \frac{\big( (\mathrm{grad} \, \kappa)(x) \big)_j \, u^{\prime}_j}{x_j}.
\end{eqnarray}
Equations~(\ref{rhs-def-grad-kappa}) and (\ref{lhs-def-grad-kappa}) 
therefore yields
\begin{eqnarray}
\label{grad-kappa-med}
\frac{\big( (\mathrm{grad} \, \kappa)(x) \big)_j}{x_j}
=
c_jx_j+{\tilde c} ,
\end{eqnarray}
where ${\tilde c}$ is a constant common in $j$ ($=1,2,\cdots,m$)
emerging from the constraint $\sum_{j=1}^{m} u^{\prime}_j=0$ to
$u^{\prime} \in T_x \regS$ (cf. the constant $c$ in Eq.~(\ref{M-G})).
Indeed, to fulfil the condition
$\sum_{j=1}^{m}\big( (\mathrm{grad} \, \kappa)(x) \big)_j =0$
for $(\mathrm{grad} \, \kappa)(x)$ to be in $T_{x}\regS$, the constant
${\tilde c}$ turns out to satisfy
\begin{eqnarray}
\label{tilde-c-med}
0= \sum_{j=1}^{m} \big( (\mathrm{grad} \, \kappa)(x) \big)_j
=
\sum_{j=1}^{m} c_jx_j^2 + {\tilde c} \sum_{j=1}^{m} x_j
=
\sum_{j=1}^{m} c_jx_j^2 + {\tilde c} ,
\end{eqnarray}
which is solved to be
\begin{eqnarray}
\label{tilde-c-fin}
{\tilde c}= - \sum_{j=1}^{m} c_jx_j^2 .
\end{eqnarray}
In the sequel, Eqs.~(\ref{grad-kappa-med}) is put together with
(\ref{tilde-c-fin}) to show (\ref{grad-kappa}).
\section*{References}
\begin{harvard}
\item[]
Amari~S and H~Nagaoka, 2000, {\it Methods of Information Geometry},
Ttanslations of Mathematical Monographs vol.191 (Providence, AMS)
Chap.~7.3.
\item[]
Bengtsson~I and {\. Z}yczkowski~K, 2006 {\it Geometry of
Quantum States} (Cambridge, Cambridge UP), p.286.
\item[]
Karmarkar~N 1984 {\it Combinatorica} {\bf 4} 373.
\item[]
Karmarkar~N 1990 {\it Mathematical Developments from Linear Programming
(Contemp. Math. {vol.114})} eds Lagarias~J~C and Todd~M~J
(Providence: AMS) 51.
\item[]
Kobayashi~S and Nomizu~K 1969 {\it Foundations of Differential
Geometry vol.2} (New York, John Wiley) p.337.
\item[]
Nakamura~Y 1992 {\it Japan J. Indust. Appl. Math.}
{\bf 9} 133.
\item[]
Nakamura~Y 1993 {\it Japan J. Indust. Appl. Math.}
{\bf 10} 179.
\item[]
Nakamura~Y 1994a {\it Japan J. Indust. Appl. Math.}
{\bf 11} 1.
\item[]
Nakamura~Y 1994b {\it Japan J. Indust. Appl. Math.}
{\bf 11} 11.
\item[]
Nielsen~M~A and Chuang~I~L 2000 {\it Quantum Computation and Qauntum
Information} (Cambridge: Cambridge University Press) Chaps 1 and 2.
\item[]
Rao~C~R and Mitra~S~K 1971 {\it Generalized Inverse of Matrices and
its Applications} (New York, John Wiley) p.6.
\item[]
Uwano~Y 2006 {\it Czech. J. Phys.} {\bf 56} 1311.
\item[]
Uwano~Y Hino~H and Ishiwatari~Y 2007 {\it Phys. Atom. Nuclei}
{\bf 70} 784.
\end{harvard}
\end{document}